\RequirePackage{fix-cm}
\documentclass[smallextended]{svjour3}       
\smartqed  
\usepackage{graphicx}
\usepackage{amsmath}
\usepackage{amsfonts}
\usepackage{amssymb}
\usepackage[colorlinks]{hyperref}
\def\chain{(X_n)_{n\geq 0}}

\def\stany{\mathcal{X}}
\def\sspace{\mathcal{X}}
\def\borel{\mathcal{B}(\stany)}

\def\1c{\mathbb{I}_C(x)}
\def\d{\textrm{d}}
\def\dx{\d x}
\def\dy{\d y}

\def\n0{n_{0}}

\def\var{{\rm Var\hskip 0.5pt}}

\def\d{{\rm d}}

\begin{document}

\title{CLTs and asymptotic variance of time-sampled Markov chains\thanks{Supported by EPSRC grants EP/G026521/1 and EP/D002060/1 and by CRiSM.}
}


\author{ Krzysztof {\L}atuszy\'{n}ski         \and
        Gareth O. Roberts 
}


\institute{K. {\L}atuszy\'{n}ski \at
              Department of Statistics\\
University of Warwick\\ CV4 7AL, Coventry, UK \\
              \email{latuch@gmail.com}           
           \and
           G. O. Roberts \at
              Department of Statistics\\
University of Warwick\\ CV4 7AL, Coventry, UK
}

\date{Received: date / Accepted: date}

\maketitle

\begin{abstract}
For a Markov transition kernel $P$ and a probability distribution $ \mu$ on nonnegative integers, a time-sampled Markov chain evolves according to the transition kernel $P_{\mu} = \sum_k \mu(k)P^k.$
In this note we obtain CLT conditions for time-sampled Markov chains and derive a spectral formula for the asymptotic variance. Using these results we compare efficiency of Barker's and Metropolis algorithms in terms of asymptotic variance.
\end{abstract}

\section{Introduction} \label{sec_intro}

Let $P$ be an ergodic transition kernel of a Markov chain $\chain$ with limiting distribution $\pi$ on $(\sspace, \borel)$ and let $f: \sspace \to \mathbb{R}$ be in $L^2(\pi).$  A typical MCMC procedure for estimating $I = \pi f:= \int_{\sspace}f(x) \pi(dx)$ would use $
\hat{I}_n := {1 \over n} \sum_{i=0}^{n-1}f(X_i).
$
Under appropriate assumptions on $P$ and $f$ a CLT holds for $\hat{I}_n,$ i.e. 
\begin{eqnarray}
\sqrt{n}(\hat{I}_n - I) & \to & \mathcal{N}(0, \sigma^2_{f,P}),
\label{eq:CLT}
\end{eqnarray}
where the constant $\sigma^2_{f,P} < \infty$ is called asymptotic variance and depends only on $f$ and $P.$

The following theorem from \cite{kipnis1986central} is a fundamental result on conditions that guarantee \eqref{eq:CLT} for reversible Markov chains. 

\begin{theorem}[\cite{kipnis1986central}] \label{thm:KV} For a reversible and ergodic Markov chain, 
and a function $f \in L^2(\pi),$ if \begin{eqnarray} Var(f,P) & := & \lim_{n \to \infty} n \var_{\pi}(\hat{I}_{n}) \;\; < \;\; \infty, \label{eq:KV} \end{eqnarray} then \eqref{eq:CLT} holds with \begin{eqnarray} \sigma^2_{f,P} & = & Var(f,P) \;\; = \;\; \int_{[-1,1]}{1+x \over 1-x} E_{f,P}(dx),
\label{eq:KV_spectral}
\end{eqnarray}
where $E_{f,P}$ is the spectral measure associated with $f$ and $P$.
\end{theorem}

We refer to \eqref{eq:KV} as the Kipnis-Varadhan condition. Assuming that \eqref{eq:KV} holds and $P$ is reversible, in Section \ref{sec:time-sampled} we obtain conditions for the CLT and derive a spectral formula for the asymptotic variance $\sigma^2_{f,P_{\mu}}$ of a time-sampled Markov chain of the form \begin{eqnarray}
P_{\mu} & := & \sum_{k=0}^{\infty} \mu(k) P^k, 
\label{eq:time_sampled}
\end{eqnarray} 
where $\mu$ is a probability distribution on the nonnegative integers. Time-sampled Markov chains are of theoretical interest in the context of petite sets (cf. Chapter 5 of \cite{meyn1993markov}), and also in the context of computational algorithms \cite{rosenthal2003asymptotic,rosenthal2003geometric}.

Next we proceed to analyze efficiency of Barker's algorithm \cite{barker1965monte}. Barker's algorithm, similarly as Metropolis, uses an irreducible transition kernel $Q$ to draw proposals. A move form $X_n = x$ to a proposal $Y_{n+1}=y$ is then accepted with probability   
\begin{eqnarray}
\alpha^{\textrm{(B)}}(x,y) & = & {\pi(y)q(y,x) \over \pi(y)q(y,x)+ \pi(x)q(x,y)},
\label{eq:Barker1}
\end{eqnarray}
where $q(x,\cdot)$ is the transition density of $Q(x, \cdot).$ It is well known that with the same proposal kernel $Q$, the Metropolis acceptance ratio results in a smaller asymptotic variance then Barker's. In Section \ref{sec:Barker} we show that the asymptotic variance of Barker's algorithm is not bigger then, roughly speaking, two times that of Metropolis. We also motivate our considerations by recent advances in exact MCMC for diffusion models. The theoretical results are illustrated by a simulation study in Section \ref{sec:sim}.

\section{Time-sampled Markov chains} \label{sec:time-sampled}

In this section we work under assumptions of Theorem \ref{thm:KV} which imply that the asymptotic variance $\sigma^2_{f,P}$ equals $Var(f,P)$ defined in \eqref{eq:KV} and satisfies \eqref{eq:KV_spectral}. For other Markov chain CLT conditions we refer to \cite{jones2004markov,roberts2004general,meyn1993markov,ECP2008-9,roberts2008variance}. 

\begin{theorem}\label{prop_time_sampled} Let $P$ be a reversible and ergodic
transition kernel with stationary measure $\pi,$ and let $f \in L^2(\pi).$ Assume that the Kipnis-Varadhan condition (\ref{eq:KV}) holds for $f$ and $P$. For a probability distribution $\mu$ on nonnegative integers, let the time-sampled kernel $P_{\mu}$ be defined by \eqref{eq:time_sampled}. Then, if any of the following conditions hold \begin{itemize}
\item[(i)] $\mu_{\textrm{odd}} := \mu(\{1,3,5,\dots\}) > 0,$ 
\item[(ii)] $\mu(0) < 1$ and $P$ is geometrically ergodic,
\end{itemize}
the CLT holds for $f$ and $P_{\mu},$ moreover \begin{eqnarray}
\label{eqn_as_var_of_time-sampled}
\sigma^2_{f, P_{\mu}} & = & \int_{[-1,1]}\frac{1+G_{\mu}(x)}{1-G_{\mu}( x)}E_{f,P}(dx)\;\; < \;\; \infty,
\end{eqnarray}
where $G_{\mu}$ is the probability generating function of $\mu,$ i.e. $G_{\mu}(z) := \mathbb{E}_{\mu}z^K,$ $|z| \leq 1, \; K \sim \mu,$ and $E_{f,P}$ is the spectral measure associated with $f$ and $P$.
\end{theorem}

\begin{remark} The condition $\mu_{\textrm{odd}} > 0$ in the above result is necessary, which we show below by means of a counterexample.
\end{remark}

\begin{proof} The proof is based on the functional analytic approach
(see e.g. \cite{kipnis1986central,roberts1997geometric}). Without loss of generality assume that $\pi f = 0.$ A reversible
transition kernel $P$ with invariant distribution $\pi$ is a
self-adjoint operator on $L_0^2(\pi):= \{ f \in L^2(\pi): \pi f =0\}$ with spectral radius bounded by
1. By the spectral decomposition theorem for self adjoint operators,
for each $f \in L^2_0(\pi)$ there exists a finite positive measure
$E_{f,P}$ on $[-1,1],$ such that \begin{eqnarray} \nonumber
\langle f, P^n f \rangle & = & \int_{[-1,1]} x^n E_{f,P}(dx),
\end{eqnarray}
for all integers $n \geq 0.$ Thus in particular \begin{eqnarray}
\label{eqn_spectral_in_proof_stat_var}
\sigma^2_f & = & \pi f^2 = \int_{[-1,1]}1 E_{f,P}(dx)\;\; <\;\; \infty, \end{eqnarray}
and by \cite{kipnis1986central} (c.f. also Theorem 4 of \cite{häggström2007variance}) one obtains 
\begin{eqnarray}
\sigma^2_{f,P} & = & \int_{[-1,1]}\frac{1+x}{1-x}E_{f,P}(dx) \;\;< \;\;\infty.
\label{eqn_spectral_in_proof_as_var}
\end{eqnarray}
Since $P_{\mu}^n \;\;= \;\; \sum_{k} \mu(k) P^k,$
by the spectral mapping theorem \cite{conway1990course}, we have
\begin{eqnarray}\nonumber
\langle f, P^n_{\mu}f \rangle & = &  \int_{[-1,1]} x^n E_{f,P_{\mu}}(dx)  \;\; = \;\; \int_{[-1,1]} \Big( \sum_k \mu(k)x^k \Big)^n E_{f,P}(dx) \\ \nonumber  & = & \int_{[-1,1]}  \Big( G_{\mu}(x)\Big)^n E_{f,P}(dx), 
\nonumber \end{eqnarray}
and consequently, applying the same argument as \cite{kipnis1986central,häggström2007variance}, we obtain
\begin{eqnarray} \nonumber
\sigma^2_{f,P_{\mu}} & = & \int_{[-1,1]}\frac{1+x}{1-x}E_{f,P_{\mu}}(dx) \\ \label{int_as_var} &= & \int_{[-1,1]}\frac{1+G_{\mu}(x)}{1-G_{\mu}(x)}E_{f,P}(dx) \label{int_for_as_var} \;\; =: \;\; \clubsuit.
\end{eqnarray}
Now \eqref{int_as_var} gives the claimed formula but we need to prove \eqref{int_as_var} is finite: by \cite{kipnis1986central} finiteness of the integral in (\ref{int_as_var}) implies a CLT for $f$ and $P_{\mu}$. Observe that \begin{eqnarray*}
|G(x)| & \leq & 1 \qquad \qquad \qquad \quad  \qquad \quad \, \;\; \textrm{for all} \quad x \in [-1,1], \\
G(x) & \leq & \mu(0) + x(1-\mu(0)) \qquad \quad  \, \textrm{for} \quad x \geq 0. \end{eqnarray*}
Moreover, if (i) holds, then \begin{eqnarray*} 
G(x) & \leq & \sum_{k \; \; \textrm{even}} \mu(k) x^k \;\; \leq \;\; 1-\mu_{\textrm{odd}} \qquad \; \; \textrm{for} \quad x \leq 0,  
\label{eq:}
\end{eqnarray*} 
hence we can write
\begin{eqnarray} \nonumber
\clubsuit & = & \int_{[-1,0)}\frac{1+G_{\mu}(x)}{1-G_{\mu}(x)}E_{f,P}(dx) + \int_{[0,1]}\frac{1+G_{\mu}(x)}{1-G_{\mu}(x)}E_{f,P}(dx) \\ 
& \leq & {1 \over \mu_{\textrm{odd}}} \int_{[-1,0)}2 E_{f,P}(dx) + {1 \over 1 - \mu(0) } \int_{[0,1]}\frac{2}{1-x}E_{f,P}(dx).
\label{eq:as_var_split_for_finiteness}
\end{eqnarray}
The first integral in \eqref{eq:as_var_split_for_finiteness} is finite by \eqref{eqn_spectral_in_proof_stat_var} and the second by \eqref{eqn_spectral_in_proof_as_var} and we are done with (i).

Next assume that (ii) holds. By $S(P)$ denote the spectrum of $P$ and let $s_P := \sup \{|\lambda|: \lambda \in S(P)\}$ be the spectral radius. From \cite{roberts1997geometric} we know that since $P$ is reversible and geometrically ergodic, it has a spectral gap, i.e. $s_P < 1.$ Hence for $x \in [-s_P, 0],$ we can write \begin{eqnarray*}
G_{\mu} & \leq & \mu(0) +  \sum_{k \; \; \textrm{even}} \mu(k)x^k \;\; \leq \;\; \mu(0) + s_P(1-\mu(0)).
\end{eqnarray*}
Consequently
\begin{eqnarray} \nonumber
\clubsuit & = & \int_{[-s_P,0)}\frac{1+G_{\mu}(x)}{1-G_{\mu}(x)}E_{f,P}(dx) + \int_{[0,s_P]}\frac{1+G_{\mu}(x)}{1-G_{\mu}(x)}E_{f,P}(dx) \\ 
& \leq & {1 \over 1-\mu(0)} \int_{[-s_P,0)}{2 \over 1- s_P} E_{f,P}(dx) + {1 \over 1 - \mu(0) } \int_{[0,s_P]}\frac{2}{1-x}E_{f,P}(dx).\qquad 
\label{eq:as_var_split_for_finiteness2}
\end{eqnarray}
The first integral in \eqref{eq:as_var_split_for_finiteness2} is finite by \eqref{eqn_spectral_in_proof_stat_var} and the second by \eqref{eqn_spectral_in_proof_as_var}.
\end{proof}


The most important special case of Theorem \ref{prop_time_sampled} is underlined and computed explicitly in the next corollary.

\begin{corollary}\label{cor_as_var_of_lazy} Let $P$ be a reversible and ergodic
transition kernel with stationary measure $\pi,$ and assume that for $f$ and $P$ the CLT (\ref{eq:CLT}) holds. For $\varepsilon \in
(0,1)$ let the lazy version of $P$ be defined as $P_{\varepsilon} \; :=
\; \varepsilon \textrm{Id} + (1-\varepsilon) P.$ Then the CLT holds for $f$ and $P_{\varepsilon}$ and \begin{eqnarray}
\label{eqn_as_var_of_lazy}
\sigma^2_{f, P_{\varepsilon}} & = & \frac{1}{1-\varepsilon}\sigma^2_{f,P} +
\frac{\varepsilon}{1-\varepsilon} \sigma^2_f.
\end{eqnarray}
\end{corollary}

\begin{proof} We use Theorem \ref{prop_time_sampled} with $\mu(0) = \varepsilon,$ $\mu(1)=1-\varepsilon.$ Hence $G_{\mu} = \varepsilon + (1-\varepsilon)x,$ and consequently 
\begin{eqnarray*}
\sigma^2_{f,P_{\varepsilon}} & = & \int_{[-1,1]}\frac{1+\varepsilon + (1-\varepsilon)
x}{1-\varepsilon - (1-\varepsilon) x}E_{f,P}(dx) \\ \nonumber
& = & \int_{[-1,1]}\frac{1}{1-\varepsilon}\bigg( \frac{1+ x}{1- x} + 
\varepsilon \bigg)E_{f,P}(dx)  \\ \nonumber & = & \frac{1}{1-\varepsilon}
\int_{[-1,1]}\frac{1+x}{1-x}E_{f,P}(dx)     +
\frac{\varepsilon}{1-\varepsilon}\int_{[-1,1]}1 E_{f,P}(dx) \\
& = & \frac{1}{1-\varepsilon}\sigma^2_{f,P} +
\frac{\varepsilon}{1-\varepsilon} \sigma^2_f.
\end{eqnarray*}
\end{proof}

Efficiency of time sampled Markov chains can be compared using the following  corollary from Theorem \ref{prop_time_sampled}.

\begin{corollary} \label{cor:eff_of_positive} Let $P$ and $f$ be as in Theorem \ref{prop_time_sampled}. If $P$ is positive as an operator on $L^2(\pi)$ and $\mu_1$ dominates stochastically $\mu_2$ (i.e. $\mu_1 \geq_{st} \mu_2$), then $P_{\mu_1}$ dominates $P_{\mu_2}$ in the efficiency ordering, i.e. $\sigma^2_{f,P_{\mu_1}} \leq \; \sigma^2_{f,P_{\mu_2}}.$
\end{corollary}
\begin{proof} If $P$ is positive self-adjoint then $\textrm{supp}E_{f,P} \subseteq [0,1].$ Moreover \begin{eqnarray*} \mu_1 \geq_{st} \mu_2 & \Rightarrow & G_{\mu_1}(x) \leq G_{\mu_2}(x) \qquad \textrm{ for} \quad x \in [-1,1].\end{eqnarray*} The conclusion follows from \eqref{eqn_as_var_of_time-sampled}.
\end{proof}

In another direction of studying CLTs, the \emph{variance bounding} property of Markov chains has been introduced in \cite{roberts2008variance} and is defined as follows. $P$ is variance bounding if there exists $K < \infty$ such that $Var(f,P) \leq K \var_{\pi}(f)$ for all $f.$ Here $Var(f,P)$ is defined in \eqref{eq:KV} and $\var_{\pi}(f) = \pi f^2 - (\pi f)^2.$ We prove that for time-sampled Markov chains the variance bounding property propagates the same way the CLT does.

\begin{theorem} Assume $P$ is reversible and variance bounding. Then $P_{\mu}$ is variance bounding if any of the following conditions hold
\begin{itemize}
\item[(i)] $\mu_{\textrm{odd}} := \mu(\{1,3,5,\dots\}) > 0,$ 
\item[(ii)] $\mu(0) < 1$ and $P$ is geometrically ergodic.
\end{itemize}
\end{theorem} 
\begin{proof} For any $f$ such that $\var_{\pi}f < \infty,$ the Kipnis-Varadhan condition holds due to variance bounding property of $P$ and thus the assumptions of Theorem~\ref{prop_time_sampled} are met. Hence for every $f \in L^2(\pi)$ there is a CLT for $f$ and $P_{\mu}.$ Therefore $P_{\mu}$ is variance bounding by Theorem 7 of \cite{roberts2008variance}.
\end{proof}

The next example shows that in case of Markov chains that are not geometrically ergodic, the condition $\mu_{\textrm{odd}} > 0$ is necessary. 

\begin{example} We set $f(x) = x$ and give an example of an ergodic and reversible transition kernel $P$ on $\mathcal{X} = [-1,1],$ and such that there is a CLT for $P$ and $f$ but \emph{not} for $P^2$ and $f.$ We shall rely on Theorem 4.1 of \cite{ECP2008-9} that provides if and only if conditions for Markov chains CLTs in terms of regenerations. It will be apparent that the condition $\mu_{\textrm{odd}}>0$ in Theorem \ref{prop_time_sampled} is necessary.

Set $s(x):= \sqrt{1-|x|},$ let $U(\cdot)$ be the uniform distribution on $[-1,1],$ and let the kernel $P$ be of the form \begin{eqnarray}
P(x, \cdot) & = & (1-s(x))\delta_{-x}(\cdot) + s(x) U(\cdot), \qquad \qquad \textrm{hence} \qquad \\
\label{eq:P_mixture}
P^2(x, \cdot) & = & (1-s(x))^2\delta_{x}(\cdot) + (2s(x)-s(x)^2) U(\cdot).
\end{eqnarray}
To find the stationary distribution of $P$ (and also $P^2$), we verify reversibility with $\pi(x) \propto 1/s(x).$
\begin{eqnarray}
\pi(\dx)P(x,\dy) & \propto & {1 \over s(x)} \delta_{-x}(y) + \delta_{-x}(y) + U(\dy) \nonumber \\
& = & {1 \over s(y)} \delta_{-y}(x) + \delta_{-y}(x) + U(\dx) \;\; \propto \;\; \pi(\dy)P(y,\dx). \nonumber 
\label{eq:}
\end{eqnarray}
Hence $\pi(x)$ is a reflected Beta$(1, \frac{1}{2}).$  Clearly $\pi(f^2) < \infty.$ 

 Recall now the split chain construction \cite{nummelin1978splitting,athreya1978new} of the bivariate Markov chain $\{X_n, \Gamma_n\}$ on $\{0,1\}\times \mathcal{X} = \{0,1\}\times [0,1].$ If $(X_n)_{n \geq 0}$ evolves according to $P$ defined in \eqref{eq:P_mixture}, we have the following transition rule from $\{X_{n-1}, \Gamma_{n-1}\}$ to  $\{X_n, \Gamma_n\}$ for the split chain.
 \begin{eqnarray*}
\check{\mathbb{P}}(X_n \in \cdot | \Gamma_{n-1} = 1, X_{n-1} = x) & = & U(\cdot), \\
\check{\mathbb{P}}(X_n \in \cdot | \Gamma_{n-1} = 0, X_{n-1} = x) & = & \delta_{-x}(\cdot), \\
\check{\mathbb{P}}(\Gamma_n =1 | \Gamma_{n-1}, X_{n} = x) & = & s(x), \\
\check{\mathbb{P}}(\Gamma_n =0 | \Gamma_{n-1}, X_{n} = x) & = & 1-s(x).
\label{eq:}
\end{eqnarray*}
The notation $\check{\mathbb{P}}$ above indicates that we consider the extended probability space for $(X_n, \Gamma_n),$ not the original one of $X_n.$ The appropriate modification of the above holds if the dynamics of $X_n$ is $P^2,$ namely 
 \begin{eqnarray*}
\check{\mathbb{P}}(X_n \in \cdot | \Gamma_{n-1} = 1, X_{n-1} = x) & = & U(\cdot), \\
\check{\mathbb{P}}(X_n \in \cdot | \Gamma_{n-1} = 0, X_{n-1} = x) & = & \delta_x(\cdot), \\
\check{\mathbb{P}}(\Gamma_n =1 | \Gamma_{n-1}, X_{n} = x) & = & 2s(x)-s^2(x), \\
\check{\mathbb{P}}(\Gamma_n =0 | \Gamma_{n-1}, X_{n} = x) & = & (1-s(x))^2.
\label{eq:}
\end{eqnarray*}
We refer to to the original papers for more details on the split chain construction and to \cite{ECP2008-9,roberts2004general} for central limit theorems in this context.  Denote \begin{eqnarray}
\tau & := & \min\{k \geq 0: \Gamma_k = 1\}.
\label{eq:tau}
\end{eqnarray} 

By Theorem 4.1 of \cite{ECP2008-9}, the CLT for $P$ and $f$ holds if and only if the following expression for the asymptotic variance is finite. 
\begin{eqnarray}
\sigma^2_{f, P}& =  & \int_{[-1,1]}s(x)\pi(x) \dx \;\check{\mathbb{E}}_{U} \Big(\sum_{k=0}^{\tau}f(X_n)\Big)^2,
\label{eq:as_var_regen_P}
\end{eqnarray}
where $(X_n, \Gamma_n)$ follow the dynamics of $P.$ Respectively, the CLT for $P^2$ and $f$ holds in our setting, \emph{if and only if}
\begin{eqnarray}
\sigma^2_{f, P^2}& =  & \int_{[-1,1]}(2s(x)-s^2(x))\pi(x) \dx \;\check{\mathbb{E}}_{U} \Big(\sum_{k=0}^{\tau}f(X_n)\Big)^2
\label{eq:as_var_regen_P2}
\end{eqnarray}
is finite, where $(X_n, \Gamma_n)$ follow the dynamics of $P^2.$ 
 
Now observe that if $(X_n)_{n \geq 0}$ evolves according to $P,$ then  $(\sum_{k=0}^{\tau}f(X_n))^2$ equals $0$ if $\tau$ is odd, or $(\sum_{k=0}^{\tau}f(X_n))^2 = X^2_0,$ if $\tau$ is even. Consequently \eqref{eq:as_var_regen_P} is finite. However, if  $(X_n)_{n \geq 0}$ evolves according to $P^2,$ then $(\sum_{k=0}^{\tau}f(X_n))^2 = (\tau+1)^2X_0^2$ and the distribution of $\tau$ is geometric with parameter $2s(X_0) - s^2(X_0) = 1 - (1-s(x))^2.$ Therefore we compute $\sigma^2_{f, P^2}$ in \eqref{eq:as_var_regen_P2} as
\begin{eqnarray}
\sigma^2_{f, P^2} & = & \int_{[-1,1]}(2s(x)-s^2(x))\pi(x) \dx \;\int_{[-1,1]}  {2 - \big(1 -\big(1-s(x)\big)^2\big)\over 2\big(1 -(1-s(x))^2\big)^2} x^2 \dx \nonumber \\
& = & C \;\int_{[-1,1]} {\big(1+(1-s(x))^2 \big) x^2 \over 2\big(1-|x| - 2\sqrt{1-|x|}\big)^2} \dx \nonumber \\
& \geq & C \;\int_{[-1,1]} { x^2 \over 8(1-|x|) } \dx \;\; = \;\; \infty. \nonumber
\label{eq:}
\end{eqnarray}

\end{example}

\section{Barker's algorithm}\label{sec:Barker}

When assessing efficiency of Markov chain Monte Carlo algorithms, the asymptotic variance criterion is one of natural choices. Peskun ordering \cite{peskun1973optimum} (see also \cite{tierney1998note,mira1999ordering}) provides a tool to compare two reversible transition kernels $P_1,$ $P_2$ with the same limiting distribution $\pi$ and is defined as follows. $P_1 \succ P_2 \iff$ for $\pi-$almost every $x \in \sspace$ and all $A \in \borel$ holds $P_1(x, A-\{x\}) \geq P_2(x, A-\{x\}).$ If  $P_1 \succ P_2$ then $\sigma^2_{f, P_1} \leq \sigma^2_{f, P_2}$ for every $f \in L^2(\pi).$

Consider now a class of algorithms where the transition kernel $P$ is defined by applying an irreducible proposal kernel $Q$ and an acceptance rule $\alpha,$ i.e. given $X_n=x,$ the value of $X_{n+1}$ is a result of performing the following two steps.
\begin{enumerate}
	\item Draw a proposal $y \sim Q(x, \cdot),$
	\item Set $X_{n+1}: = y$ with probability $\alpha(x,y)$ and $X_{n+1} = x$ otherwise, 
\end{enumerate}
where $\alpha(x,y)$ is such that the resulting kernel $P$ is reversible with stationary distribution $\pi$. It follows \cite{peskun1973optimum,tierney1998note} that for a given proposal kernel $Q$ the standard Metropolis-Hastings \cite{metropolis1953equations,hastings1970monte} acceptance rule \begin{eqnarray}
\alpha^{\textrm{(MH)}}(x,y) & = & \min\big\{1, {\pi(y)q(y,x) \over \pi(x)q(x,y)} \big\}
\label{eq:MH}
\end{eqnarray}
yields a transition kernel $P^{\textrm{(MH)}}$ that is maximal with respect to Peskun ordering and thus minimal with respect to asymptotic variance. In particular, the Barker's algorithm \cite{barker1965monte} that uses acceptance rule 
\begin{eqnarray}
\alpha^{\textrm{(B)}}(x,y) & = & {\pi(y)q(y,x) \over \pi(y)q(y,x)+ \pi(x)q(x,y)}
\label{eq:Barker's}
\end{eqnarray}
is inferior to Metropolis-Hastings when the asymptotic variance is considered.
In the above notation we assume that all the involved distributions have common denominating measure and $q(x,\cdot)$ are transition densities of $Q.$ See \cite{tierney1998note} for a more general statement and discussion.

Exact Algorithms introduced in \cite{beskos2006exact,beskos2005exact,beskos2006retrospective,beskos2008factorisation} allow for inference in diffusion models without Euler discretization error. In recent advances in Exact MCMC inference for complex diffusion models a particular setting is reoccurring, where the Metropolis-Hastings acceptance step requires a specific Bernoulli Factory and is not possible to execute. However, in this diffusion context the Barker's algorithm (\ref{eq:Barker's}) is feasible, as well as the 'lazy' version of the Metropolis-Hastings kernel \begin{eqnarray}
P_{\varepsilon}^{\textrm{(MH)}} &:= &\varepsilon Id + (1-\varepsilon)P^{\textrm{(MH)}}.\label{eq:lazy}\end{eqnarray}
We refer to \cite{Flavio_jump,LaPaRo_MarkovSwitching,p2p} for the background on exact MCMC inference for diffusions and the Bernoulli Factory problem. This motivates us to investigate performance of these alternatives in comparison to the standard Metropolis-Hastings.

\begin{theorem}\label{thm_as_var} Let $P^{\textrm{(B)}}$ denote the transition kernel of the Barker's algorithm and let $P^{\textrm{(MH)}}$ and $P^{\textrm{(MH)}}_{\varepsilon}$ be as defined in (\ref{eq:lazy}). If the CLT (\ref{eq:CLT}) holds for $f$ and $P^{\textrm{(MH)}},$ then it holds also for \begin{itemize}
   \item[(i)] $f$ and $P^{\textrm{(MH)}}_{\varepsilon}$ with 
   \begin{eqnarray}
\label{eqn_as_var_MH}
\sigma^2_{f, P^{\textrm{(MH)}}_{\varepsilon}} & = & \frac{1}{1-\varepsilon}\sigma^2_{f,P^{\textrm{(MH)}}} +
\frac{\varepsilon}{1-\varepsilon} \sigma^2_f.
\end{eqnarray}
   \item[(ii)] $f$ and $P^{\textrm{(B)}}$ with 
   \begin{eqnarray}
\label{eqn_as_var_B}
\sigma^2_{f,P^{\textrm{(MH)}}} &\leq & \sigma^2_{f, P^{\textrm{(B)}}} \;\; \leq \;\; \sigma^2_{f, P^{\textrm{(MH)}}_{1/2}} \;\; = \;\;  2 \sigma^2_{f,P^{\textrm{(MH)}}} +
 \sigma^2_f.
\end{eqnarray}
\end{itemize}
\end{theorem}
 
 \begin{proof}
  The first claim $(i)$ is a restatement of Corollary \ref{cor_as_var_of_lazy} for Metropolis-Hastings chains. To obtain the second claim $(ii),$ note that $P^{\textrm{(MH)}}_{1/2}$ can be viewed as an algorithm that uses proposals from $Q$ and acceptance rule \begin{eqnarray}
\alpha(x,y) & = & \min\big\{ {1\over 2}, {\pi(y)q(y,x) \over 2 \pi(x)q(x,y)} \big\}.
\nonumber
\end{eqnarray} Now since \begin{eqnarray}
\min\big\{ 1, {\pi(y)q(y,x) \over \pi(x)q(x,y)} \big\} & \geq & {\pi(y)q(y,x) \over \pi(y)q(y,x) + \pi(x)q(x,y) }   \; \geq \; \min\big\{ {1\over 2}, {\pi(y)q(y,x) \over 2 \pi(x)q(x,y) } \big\}, \nonumber
\label{eq:}
\end{eqnarray} the result follows from Peskun ordering and Corollary \ref{cor_as_var_of_lazy}.
 \end{proof}

\section{Numerical Examples} \label{sec:sim}

To illustrate the theoretical findings, we consider two numerical examples. The first focuses on time sampling, the second on efficiency of the Barker's algorithm.

\subsection{Time sampled contracting normals} 
Consider the contracting normals example, i.e. a Markov chain with transition probabilities \begin{equation}
P(x, \cdot) = N(\theta x , 1- \theta^2)
\label{eq:}
\end{equation}
for some $\theta \in (-1,1).$
It is easy to check that the stationary distribution is $\pi(\cdot) = N(0,1).$ Moreover the transition kernel is geometrically ergodic and reversible for all $\theta \in(-1,1)$ and also positive for $\theta \in [0,1),$ \cite{baxendale2005renewal,eps_alpha}. For the target function we take $f(x) =x$ and estimate the asymptotic variance using the batch means estimator of \cite{jones2006fixed} based on a trajectories of length $10^7.$ We set $\theta$ to  $0.9$ and $-0.9$ in the following settings: 
\begin{itemize}
	\item CN: Contracting normals; 
	\item LCN: Lazy contracting normals with 
						$\varepsilon = 0.5$;
	\item TSCN1: Time sampled contracting normals for 
						sampling distribution $$\mu = 1+ Pois(1);$$
	\item TSCN2: Time sampled contracting normals for 
						sampling distribution $$\mu = 1+ Pois(5).$$
\end{itemize}
\begin{center}
\begin{tabular}{|c|c|c|c|c|}
\hline
 & CN &  LCN & TSCN1 &TSCN2  \\
\hline
$\theta = 0.9$  & 19.1 & 38.5 & 9.28   & 3.43  \\
\hline
$\theta = - 0.9$ & 0.053 & 1.14 & 0.80  & 0.96  \\
\hline
\end{tabular} 
{\footnotesize
\smallskip

 Table 1. Estimated asymptotic variance of the contracting normals Markov chain for different sampling scenarios. 
}
\end{center}

The first two columns of Table 1 report how laziness increases asymptotic variance and illustrate Corollary \ref{cor_as_var_of_lazy}. Note that the stationary variance $\sigma^2_f = 1$ is substantial compared to the asymptotic variance of contracting normals for $\theta = -0.9$ and thus the lazy version LCN becomes severely inefficient compared to CN. The stochastic ordering of the sampling distributions  in the above scenarios is LCN $<_{st}$ CN $<_{st}$ TSCN1 $<_{st}$ TSCN2 therefore the simulation shows how the asymptotic variance decreases for stochastically bigger sampling distributions (Corollary \ref{cor:eff_of_positive}) in case of positive operators ($\theta = 0.9$) and how this property fails if the operator is not positive, i.e for $\theta = -0.9.$

\subsection{Efficiency of the Barker's algorithm}

We compare the estimated asymptotic variance of the random walk Metropolis algorithm, the Barker's algorithm and lazy version of the random walk Metropolis with $\varepsilon = 0.5$ to illustrate the bounds of Theorem \ref{thm_as_var}. For the stationary distribution we take $N(0,1)$ and the increment proposal is $U([-2,2]).$ The results based on a simulation length $10^7$ are reported in Table 2.
\begin{center}
\begin{tabular}{|c|c|c|c|}
\hline
 & Metropolis &  Barker's & lazy Metropolis   \\
\hline
asymptotic variance  & 3.69  & 5.67   & 8.32  \\
\hline
\end{tabular} 
{\footnotesize
\smallskip

 Table 1. Estimated asymptotic variance of the Metropolis, Barker's and lazy Metropolis algorithms. 
}
\end{center}

\section{Acknowledgements}
We thank Jeffrey S. Rosenthal for a helpful discussion.

\bibliographystyle{spmpsci} 
\bibliography{references}

\begin{thebibliography}{10}
\providecommand{\url}[1]{{#1}}
\providecommand{\urlprefix}{URL }
\expandafter\ifx\csname urlstyle\endcsname\relax
  \providecommand{\doi}[1]{DOI~\discretionary{}{}{}#1}\else
  \providecommand{\doi}{DOI~\discretionary{}{}{}\begingroup
  \urlstyle{rm}\Url}\fi

\bibitem{athreya1978new}
Athreya, K., Ney, P.: {A new approach to the limit theory of recurrent Markov
  chains}.
\newblock Transactions of the American Mathematical Society \textbf{245}(Nov),
  493--501 (1978)

\bibitem{barker1965monte}
Barker, A.: {Monte Carlo calculations of the radial distribution functions for
  a proton-electron plasma}.
\newblock Australian Journal of Physics \textbf{18}, 119 (1965)

\bibitem{baxendale2005renewal}
Baxendale, P.: {Renewal theory and computable convergence rates for
  geometrically ergodic Markov chains}.
\newblock The Annals of Applied Probability \textbf{15}(1B), 700--738 (2005)

\bibitem{ECP2008-9}
Bednorz, W., {\L}atuszy\'nski, K., Lata{\l}a, R.: A regeneration proof of the
  central limit theorem for uniformly ergodic markov chains.
\newblock Electronic Communications in Probability \textbf{13}, 85--98 (2008)

\bibitem{beskos2006retrospective}
Beskos, A., Papaspiliopoulos, O., Roberts, G.: {Retrospective exact simulation
  of diffusion sample paths with applications}.
\newblock Bernoulli \textbf{12}(6), 1077 (2006)

\bibitem{beskos2008factorisation}
Beskos, A., Papaspiliopoulos, O., Roberts, G.: {A factorisation of diffusion
  measure and finite sample path constructions}.
\newblock Methodology and Computing in Applied Probability \textbf{10}(1),
  85--104 (2008)

\bibitem{beskos2006exact}
Beskos, A., Papaspiliopoulos, O., Roberts, G., Fearnhead, P.: {Exact and
  computationally efficient likelihood-based estimation for discretely observed
  diffusion processes (with discussion)}.
\newblock Journal of the Royal Statistical Society: Series B (Statistical
  Methodology) \textbf{68}(3), 333--382 (2006)

\bibitem{beskos2005exact}
Beskos, A., Roberts, G.: {Exact simulation of diffusions}.
\newblock Annals of Applied Probability \textbf{15}(4), 2422--2444 (2005)

\bibitem{conway1990course}
Conway, J.: {A course in functional analysis}.
\newblock Springer (1990)

\bibitem{Flavio_jump}
Gon\c{c}alves, F., Roberts, G., {\L}atuszy\'{n}ski, K.: Exact mcmc inference
  for jump diffusion models with stochastic jump rate (2011)

\bibitem{häggström2007variance}
H{\" a}ggstr{\" o}m, O., Rosenthal, J.: {On variance conditions for Markov
  chain CLTs}.
\newblock Elect. Comm. in Probab \textbf{12}, 454--464 (2007)

\bibitem{hastings1970monte}
Hastings, W.: {Monte Carlo sampling methods using Markov chains and their
  applications}.
\newblock Biometrika \textbf{57}(1), 97 (1970)

\bibitem{jones2004markov}
Jones, G.: {On the Markov chain central limit theorem}.
\newblock Probability surveys \textbf{1}, 299--320 (2004)

\bibitem{jones2006fixed}
Jones, G., Haran, M., Caffo, B., Neath, R.: {Fixed-width output analysis for
  Markov chain Monte Carlo}.
\newblock Journal of the American Statistical Association \textbf{101}(476),
  1537--1547 (2006)

\bibitem{kipnis1986central}
Kipnis, C., Varadhan, S.: {Central limit theorem for additive functionals of
  reversible Markov processes and applications to simple exclusions}.
\newblock Communications in Mathematical Physics \textbf{104}(1), 1--19 (1986)

\bibitem{p2p}
{\L}atuszy{\'n}ski, K., Kosmidis, I., Papaspiliopoulos, O., Roberts, G.:
  {Simulating events of unknown probabilities via reverse time martingales}.
\newblock Random Structures \& Algorithms  (2011)

\bibitem{eps_alpha}
{\L}atuszy{\'n}ski, K., Niemiro, W.: {Rigorous confidence bounds for MCMC under
  a geometric drift condition}.
\newblock Journal of Complexity \textbf{27}(1), 23--38 (2011)

\bibitem{LaPaRo_MarkovSwitching}
{\L}atuszy\'{n}ski, K., Palczewski, J., Roberts, G.: Exact inference for a
  markov switching diffusion model with discretely observed data (2011)

\bibitem{metropolis1953equations}
Metropolis, N., Rosenbluth, A., Rosenbluth, M., Teller, A., Teller, E.:
  {Equations of state calculations by fast computational machine}.
\newblock Journal of Chemical Physics \textbf{21}(6), 1087--1091 (1953)

\bibitem{meyn1993markov}
Meyn, S., Tweedie, R.: {Markov chains and stochastic stability}.
\newblock Springer London et al. (1993)

\bibitem{mira1999ordering}
Mira, A., Geyer, C.: {Ordering Monte Carlo Markov chains}.
\newblock In: School of Statistics, University of Minnesota. technical report
  (1999)

\bibitem{nummelin1978splitting}
Nummelin, E.: {A splitting technique for Harris recurrent Markov chains}.
\newblock Probability Theory and Related Fields \textbf{43}(4), 309--318 (1978)

\bibitem{peskun1973optimum}
Peskun, P.: {Optimum monte-carlo sampling using markov chains}.
\newblock Biometrika \textbf{60}(3), 607 (1973)

\bibitem{roberts1997geometric}
Roberts, G., Rosenthal, J.: {Geometric ergodicity and hybrid Markov chains}.
\newblock Electron. Comm. Probab \textbf{2}(2), 13--25 (1997)

\bibitem{roberts2004general}
Roberts, G., Rosenthal, J.: {General state space Markov chains and MCMC
  algorithms}.
\newblock Probability Surveys \textbf{1}, 20--71 (2004)

\bibitem{roberts2008variance}
Roberts, G., Rosenthal, J.: {Variance bounding Markov chains}.
\newblock Annals of applied probability \textbf{18}(3), 1201 (2008)

\bibitem{rosenthal2003asymptotic}
Rosenthal, J.: {Asymptotic variance and convergence rates of nearly-periodic
  Markov chain Monte Carlo algorithms}.
\newblock Journal of the American Statistical Association \textbf{98}(461),
  169--177 (2003)

\bibitem{rosenthal2003geometric}
Rosenthal, J.: {Geometric convergence rates for time-sampled markov chains}.
\newblock Journal of Theoretical Probability \textbf{16}(3), 671--688 (2003)

\bibitem{tierney1998note}
Tierney, L.: {A note on Metropolis-Hastings kernels for general state spaces}.
\newblock Annals of Applied Probability \textbf{8}(1), 1--9 (1998)

\end{thebibliography}

\end{document}